\documentclass[12pt]{amsart}
\usepackage{amsmath,amsthm,amsfonts,amssymb,latexsym,enumerate,xcolor}
\usepackage[pagebackref]{hyperref}
\usepackage[english]{babel}

\newcommand{\N}{\mathbb{N}}

\newcommand{\Aut}{{{\operatorname{Aut}}}}

\newcommand{\Irr}{{{\operatorname{Irr}}}}

\newcommand{\Syl}{\operatorname{Syl}}

\newcommand{\T}{\operatorname{T}}
\newcommand{\acd}{\operatorname{acd}}

\newcommand{\Pro}{\operatorname{Pr}}

\newtheorem{thm}{Theorem}[section]
\newtheorem{lem}[thm]{Lemma}
\newtheorem{con}[thm]{Conjecture}

\newtheorem{cor}[thm]{Corollary}

\newtheorem*{thmA}{Theorem A}
\newtheorem*{conA'}{Conjecture A'}
\newtheorem*{thmB}{Theorem B}
\newtheorem*{thmC}{Theorem C}
\newtheorem*{thmD}{Theorem D}

\theoremstyle{definition}
\newtheorem{rem}[thm]{Remark}

\newtheorem{note}[thm]{Note}

\numberwithin{equation}{section}

\begin{document}

\title[Supersolvability, nilpotency, commuting  probability and  average character degree]
{Supersolvability and nilpotency in terms of the commuting  probability and the average character degree}

\author{Juan Mart\'{\i}nez}
\address{Departament de Matem\`atiques, Universitat de Val\`encia, 46100
  Burjassot, Val\`encia, Spain}
\email{Juan.Martinez-Madrid@uv.es}

\thanks{Research supported by Generalitat Valenciana CIAICO/2021/163 and CIACIF/2021/228. }




\begin{abstract}
Let $p$ be a prime and let $G$ be a finite group such that the smallest prime that divides $|G|$ is $p$.  We find sharp  bounds, depending on $p$, for the commuting probability and the average character degree to guarantee that $G$ is nilpotent or supersolvable.
\end{abstract}

\maketitle



\section{Introduction}
Let $G$ be a finite group. We define the commuting probability of $G$ as the probability that two random elements of $G$ commute. That is
$$\Pro(G)=\frac{\{(x,y) \in G \times G| xy=yx\}}{|G|^2}.$$

Gustafson  \cite{Gustafson} proved that 
$$\Pro(G)=\frac{k(G)}{|G|},$$
where $k(G)$ is the number of conjugacy classes of $G$. There are many results that prove that the structure of $G$ is more restricted as $\Pro(G)$ increases. The first one, is the classical result generally attributed to Gustafson \cite{Gustafson}, which states  that if $\Pro(G)>\frac{5}{8}$, then $G$ is abelian. However, a more general result was already known. More precisely, it was proved in  \cite{J} that  if $p$ is the smallest prime dividing $|G|$ and $\Pro(G)>\frac{p^2+p-1}{p^3}$, then $G$ is abelian (see \cite[Theorem 5.2]{Orig} for an updated proof of this result).

On other side, Lescot \cite{Le} proved that if $\Pro(G)>\frac{1}{2}$, then $G$ is nilpotent. Moreover, Guralnick and Robinson \cite{GR} observed that  if $\Pro(G)>\frac{1}{p}$, where $p$ is the smallest prime dividing $|G|$, then $G$ is nilpotent.  Since $\Pro(\mathsf{S}_3)=\frac{1}{2}$, the bound $\frac{1}{2}$ cannot be improved when $p=2$, that is when $|G|$ is even. However, for $p>2$  the bound $\Pro(G)>\frac{1}{p}$ was not sharp. For example, $\mathsf{C}_7\rtimes \mathsf{C}_3$ is a non-nilpotent group with $\Pro(\mathsf{C}_7\rtimes \mathsf{C}_3)=\frac{5}{21}<\frac{1}{3}$.

In addition, Barry, MacHale and Ni Shé \cite{BMNS} proved that if $\Pro(G)>\frac{1}{3}$, then $G$ is supersolvable. We have to remark that this result was extended and the proof was simplified  by Lescot, Hung and Yang \cite{LNY}.

Our first goal in this paper is to determine the best possible functions $g_{n}(p)$ and $g_{s}(p)$ such that if $\Pro(G)>g_{n}(p)$, where $p$ is the smallest prime dividing $|G|$, then $G$ is nilpotent and if $\Pro(G)>g_{s}(p)$, then $G$ is supersolvable.  The definition of these functions depends on some concepts in number theory and will be postponed to Section \ref{com}.

\begin{thmA}
Let $G$ be a group such that $p>2$ is the smallest prime dividing $|G|$. If  $\Pro(G)> g_{n}(p)$, then $G$ is nilpotent. Moreover, the bound is sharp.
\end{thmA}

\begin{thmB}
Let $G$ be a group such that $p>2$ is the smallest prime dividing $|G|$. If   $\Pro(G)> g_{s}(p)$, then $G$ is supersolvable. Moreover, the bound is sharp.
\end{thmB}

Following \cite{ILM} given  a finite group $G$, we define the average character degree of $G$ as
$$\acd(G)=\frac{\sum_{\chi \in \Irr(G)}\chi(1)}{|\Irr(G)|},$$
where $\Irr(G)$ denotes the set of irreducible complex characters of $G$. There exists a well-known relationship between $\acd(G)$ and $\Pro(G)$
$$\frac{1}{\Pro(G)}=\frac{|G|}{|\Irr(G)|}=\frac{\sum_{\chi \in \Irr(G)}\chi(1)^2}{|\Irr(G)|}\geq \frac{(1/|\Irr(G)|)(\sum_{\chi \in \Irr(G)}\chi(1))^2}{|\Irr(G)|}=\acd(G)^2,$$
where the inequality follows from  Cauchy-Schwarz inequality. As in the case of $\Pro(G)$, there are many results that prove that  the structure of $G$ is more restricted as $\acd(G)$ decreases. Isaacs, Loukaki and Moretó \cite{ILM} proved that if $\acd(G)<\frac{3}{2}$, then $G$ is supersolvable and if $\acd(G)<\frac{4}{3}$, then $G$ is nilpotent.   On the other hand, Moretó and Hung \cite{M} proved that if $\acd(G)<\frac{16}{5}$, then $G$ is solvable (this was conjectured in \cite{ILM}).   Since $\acd(\mathsf{S}_3)=\frac{4}{3}$, $\acd(\mathsf{A}_4)=\frac{3}{2}$ and $\acd(\mathsf{A}_5)=\frac{16}{5}$, we have that these bounds cannot be improved.

Our second goal in this paper is to give an analogous version of Theorems A and B for the average character degree. More precisely, our goal is to improve on the following result.

\begin{thm}[Theorem D of \cite{ILM}]\label{original}
Let $G$ be a finite group and let $p$ be the smallest prime divisor of $|G|$. Then
\begin{itemize}

\item[a)] If $p>2$ and $\acd(G)<\frac{27}{11}$, then $G$ is supersolvable.

\item[b)] If $\acd(G)<\frac{3p}{p+2}$, then $G$ is nilpotent.
\end{itemize}
\end{thm}

First, we observe that  the bound given by a) of Theorem \ref{original}  is attained in  the group $(\mathsf{C}_5 \times \mathsf{C}_5)\rtimes \mathsf{C}_3$. However, this is a global bound, which does not depend on the smallest prime dividing $|G|$. Similarly, we also observe that if $2p+1$ is a prime (that is, if $p$ is a Sophie Germain prime), then the bound given by  b) of Theorem \ref{original}  cannot be improved since $\acd(\mathsf{C}_{2p+1}\rtimes \mathsf{C}_p)=\frac{3p}{p+2}$. However, this bound it is not best possible for all primes.  Thus, our second goal in this paper is to determine   the best possible functions $h_{n}(p)$ and $h_{s}(p)$ such that if $\Pro(G)>g_{n}(p)$, where $p$ is the smallest prime dividing $|G|$, then $G$ is nilpotent and if $\Pro(G)>g_{s}(p)$, then $G$ is supersolvable. As in the case of our results on the commuting probability, the definition of our sharp bounds depends on non-trivial results in number theory and will be postponed to Section \ref{acd}.


\begin{thmC}
Let $G$ be finite group and let $p>2$ be the smallest prime dividing $|G|$. If $\acd(G)<h_{n}(p)$, then $G$ is nilpotent.  Moreover, the bound is sharp.
\end{thmC}

\begin{thmD}
Let $G$ be finite group and let $p>3$ be the smallest prime dividing $|G|$. If $\acd(G)<h_{s}(p)$, then $G$ is supersolvable. Moreover, the bound is sharp.
\end{thmD}


The paper is organized as follows. In Section \ref{preliminary} we recall and prove some preliminary results about $\Pro(G)$ and $\acd(G)$. In Subections \ref{nilpotent} and  \ref{supersolvable} we give the proofs of Theorems A and B, respectively. In Section \ref{acd} we give the proof of Theorems C and D. Finally, in Section \ref{numbers}, we will discuss  the number-theoretical questions that arise in this paper.

\section{Preliminary results}\label{preliminary}
In this section we recall some results that we will use in the remaining (sometimes without citing them explicitly). We have separated the results on the commuting probability and the results on the average character degree

\subsection{Commuting probability}

Our first result relates the commuting probability of a group with the commuting probability of its subgroups and quotients.

\begin{lem}[(i) and (ii) of Lemma 2 of \cite{GR}]\label{basic}
Let $G$ be a finite group. Then

\begin{itemize}
\item[a)] If $H \leq G$, then $\Pro(G)\leq \Pro(H)$.

\item[b)] If $N \trianglelefteq G$, then $\Pro(G)\leq \Pro(N)\Pro(G/N)\leq \Pro(G/N)$.
\end{itemize}

\end{lem}

 We will use  the following lemma to bound $|G'|$ in terms of $\Pro(G)$.

\begin{lem}[Lemma 2.2 of  \cite{HMM}]\label{bound}
	If $p$ is the smallest prime dividing the order of a finite group $G$, then $$\Pro(G) \leq
	\frac{1 + (p^{2} -1)/|G'|}{p^2}.$$
	Moreover, equality holds if and only if all non-linear irreducible characters have degree $p$.
	\end{lem}

\begin{rem}
By using that $\Pro(G)=k(G)/|G|$ it is easy to see that if all non-linear characters have degree $f>1$, then $\Pro(G) =\frac{1 + (f^{2} -1)/|G'|}{f^2}$.
\end{rem}

\subsection{Average character degree}

We begin the results on $\acd(G)$ by proving the following easy lemma.

\begin{lem}
Let $G,H$ be  finite groups. Then $\acd(H\times G)=\acd(H)\acd(G)$.
\begin{proof}
It follows trivially using that $\Irr(H\times G)=\{\chi\times \psi|\chi \in \Irr(H), \psi \in \Irr(G)\}$. 
\end{proof}
\end{lem}

Let us introduce some notation. Let $N \trianglelefteq G$ and let $\theta \in \Irr(N)$. As usual, $\Irr(G|\theta)$ will denote the set of irreducible characters of $G$ lying over $\theta$. We set
$$\acd(G|\theta)=\frac{\sum_{\chi \in \Irr(G|\theta)}\chi(1)}{|\Irr(G|\theta)|}.$$

The following result was communicated to me by A. Moret\'o.

\begin{lem}[Moret\'o]\label{dec}
Let $G$ be a group and let $p$ be the smallest prime dividing $|G|$. If $\acd(G)\leq p$, then $\acd(G/N)\leq \acd(G)$ for every $N\trianglelefteq G$.
\begin{proof}
Suppose that  $\theta \in \Irr(N)$ is either not linear or not extendible, in both cases,  we get that $\chi(1)\geq p$ for every $\chi \in \Irr(G|\theta)$ by Clifford's correspondence (see  \cite[Theorem 6.11 ]{Isaacscar}). Therefore, $\acd(G|\theta)\geq p \geq \acd(G)$.

It follows that there exists $\lambda \in \Irr(N)$ such that $\lambda$ is linear and extendible to $\mu\in  \Irr(G)$ and $\acd(G|\lambda)\leq \acd(G)$. By Gallagher's Theorem (see \cite[Corollary 6.17 ]{Isaacscar}), we have that $\Irr(G|\lambda)=\{\mu \rho| \rho \in \Irr(G/N)\}$ and hence $\acd(G|\lambda)=\acd(G/N)$. Thus, if we assume that $\acd(G/N)>\acd(G)$, then we have that $\acd(G|\lambda)=\acd(G/N)>\acd(G)$, which is impossible. Therefore, $\acd(G/N)\leq\acd(G)$ and hence, the result follows.
\end{proof}
\end{lem}

Note that Lemma \ref{dec} provides a partial positive answer to   \cite[Question 5.1]{M2}, which  asks if the inequality $\acd(G/N)\leq \acd(G)$ holds for every group $G$ and every normal subgroup $N$ of $G$. We will also need the following result from \cite{ILM}.

\begin{thm}[Theorem 3.2 of \cite{ILM}]\label{thm32}
Let $A\triangleleft G$, where $A$ is abelian and $G$ splits over $A$. If $r$ is the number of orbits of the action of $G$ on $\Irr(A)\setminus \{1_A\}$, then there exists $t$, the size of one of these orbits such that $\acd(G)\geq \frac{t(r+1)}{t+r}$.
\end{thm}

Let us define the function 
$$f(t,r)=\frac{t(r+1)}{t+r}.$$
We observe that for $t,r\geq 0$, $f(t,r)$ is increasing in each of the variables. This is, if $t\geq t_0$ (respectively, if $r\geq r_0$), then $f(t,r)\geq f(t_0,r)$ for every $r$ (resp.  $f(t,r)\geq f(t,r_0)$ for every $t$).


We finish this section by proving a lemma which was observed in the proof of \cite[Theorem D]{ILM}. To prove it, we will use the fact that a group of odd  order does not posses any non-principal real character. This fact was proved by W. Burnside and can be found in \cite[Problem 3.16]{Isaacscar}.

\begin{lem}[Isaacs, Loukaki, Moret\'o]\label{reven}
Let $G$ be a group of odd order,  let $A\leq G$ be normal and abelian and let $r$ be the number of orbits of the action of $G$ on $\Irr(A)\setminus \{1_A\}$. Then $r$ is even and if moreover $r=2$, then both orbits have the same size.
\begin{proof}
We observe that complex conjugation permutes the orbits of the action. We claim that complex conjugation cannot fix any orbit. Since $|A|$ is odd, any of the characters in $\Irr(A)\setminus \{1_A\}$ is real and hence complex conjugation cannot fix any non-principal character of $A$. Thus, the size of an orbit fixed by the conjugation must be even, which is impossible. Thus, complex conjugation cannot fix any orbit and hence $r$ is even. If $r=2$, then the two orbits are conjugate and hence they have the same size.
\end{proof}
\end{lem}

\section{Commuting probability}\label{com}

We start by defining the functions $g_n(p)$ and $g_s(p)$ that appear in the statements of Theorem A and B.
 Let $p,q$ be two primes and $r\geq1$. We say that $p$ is a Zsigmondy's prime for $\langle q,r \rangle$ if $p$ divides $q^r-1$ but does not divide $q^k-1$ for $k<r$.   Given  a prime $p>2$, and $l \geq 1$ we define
$$\T(p,l)=\{q^r| q>p, r \geq l, q \text{ is a prime, } p  \text{ is a Zsigmondy's prime for } \langle q,r \rangle\}.$$

   With this notation, we set $t(p)=\min \T(p,1)$ and $r(p)=\min \T(p,2)$. We recall that Dirichlet's Theorem (see Theorem 7.9 of \cite{A}) asserts that given $a,n \in \N$ with $(a,n)=1$, then there exist infinitely many primes of the form $k\cdot n+a$ with $k \in \N$. As a consequence, we have that there exists a prime $q$ such that $p$ divides $q+1$. Since $p>2$, we deduce that $p$ cannot divide $q-1$ and hence $q \in \T(p,2)\subseteq \T(p,1)$. Therefore,   $\T(p,2)$ and $\T(p,1)$ are both non-empty and hence both $t(p)$ and $r(p)$ are well defined.  Now, we define some functions depending on $p$
$$f_{n}(p)=\frac{1+\frac{p^2-1}{t(p)}}{p^2}$$
and
$$f_{s}(p)=\frac{1+\frac{p^2-1}{r(p)}}{p^2}.$$

Since $t(p),r(p)\geq p+1$, we have that both functions are bounded above by $\frac{1}{p}$. Therefore, both $f_{n}(p)$ and $f_{s}(p)$ tend to $0$ as $p$ tends to infinity and thus, we can  define 
$$g_{n}(p)=\max_{q\geq p, q \text{ is prime}}\{f_{n}(q)\}$$
and
$$g_{s}(p)=\max_{q\geq p, q \text{ is prime}}\{f_{s}(q)\}.$$

We had to introduce the maximum of all $f_{n}(p)$ and $f_{s}(p)$ because the sequences $\{f_{n}(p)\}$ and  $\{f_{s}(p)\}$ are not decreasing. As an example, we have that $f_{n}(19)=\frac{29}{3629}<\frac{25}{1081}=f_{n}(23)$ and $f_{s}(29)=\frac{1061}{867941}<\frac{151}{115351}=f_{n}(31)$.

\subsection{Nilpotent odd groups}\label{nilpotent}

Our goal in this subsection is to prove Theorem A.  We will assume that  all groups in this subsection will have odd order and, in particular, they will be solvable. We recall that  $G$  is a minimal non-nilpotent group if $G$ is non-nilpotent but every proper subgroup is. We have the following result about these groups.

\begin{thm}[Theorem 9.1.9 of \cite{R}]\label{OJ}
Let $G$ be a minimal non-nilpotent group. Then $G=Q\rtimes P$, where $P\in \Syl_p(G)$ is cyclic and $Q \in \Syl_q(G)$.
\end{thm}

Now, we reduce the problem to bound $\Pro(T)$ for some groups $T$.  Let $n$ be a positive  integer. Through the rest of the paper, we will write $\mathsf{C}_{n}$ to denote the  cyclic group of order $n$.  If moreover, $n$ is the power of a prime, we will write $\mathsf{H}_{n}$ to denote the elementary abelian group of order $n$. 

\begin{lem}\label{red2}
Assume that $G$ is a non-nilpotent group such that  the smallest prime dividing $|G|$ is at least $p$, $\Pro(G)>g_{n}(p)$ and it has minimal order with these properties. Then $G$ has the form $G=\mathsf{H}_{q^l}\rtimes \mathsf{C}_r$ such that $q,r\geq p$ are two odd primes, $l \geq 1$  and the action of $\mathsf{C}_r$ is faithful and simple (that is, $\mathsf{C}_r$ does not fix any proper subspace of $\mathsf{H}_{q^l}$).
\begin{proof}
By minimality of $|G|$ and Lemma \ref{basic}, we deduce that every proper subgroup of $G$ is nilpotent and every proper quotient is nilpotent. Since all proper subgroups are nilpotent, we have that $G$ is a minimal non-nilpotent group. Thus, by Theorem \ref{OJ}, we have that $G=Q\rtimes R$, where  $R=\mathsf{C}_{r^k}$ for some $k \geq 1$.

Now, if $Z(G)>1$, then $G/Z(G)$ is nilpotent, which is impossible. Thus, we deduce that $Z(G)=1$. If $k \geq 2$, then  $\mathsf{C}_{r^{k-1}}Q$ is nilpotent and hence $1\not=\mathsf{C}_{r^{k-1}}$ is central in $G$. Thus, we deduce that $R=\mathsf{C}_r$.

Now, it only remains to prove that $Q$ is a minimal normal subgroup of $G$. Assume that there exists $1<N< Q$ a minimal normal subgroup of $G$. It follows that $G/N$ is nilpotent and hence the action of $R$ on $Q/N$  is trivial. Thus, by Corollary 3.28 of \cite{Isaacs}, we have that
$$Q/N=C_{Q/N}(R)=C_{Q}(R)/N,$$
and hence $C_{Q}(R)=Q$, or equivalently  $[Q,R]=1$. Since $G=Q\rtimes R$, we deduce that $G$ is nilpotent, which is impossible.  That contradiction implies that $Q$ is a minimal normal subgroup and the result follows.
\end{proof}
\end{lem}

Now, we only have to bound $\Pro(T)$ for groups $T$ of the form described in Lemma \ref{red2}. We define $T_p:=\mathsf{H}_{t(p)}\rtimes \mathsf{C}_p$. We observe that $T_p$ has the form described in Lemma \ref{red2} and $\Pro(T_p)=f_{n}(p)$. We have  the following result.

\begin{lem}\label{prev2}
 Let $G= \mathsf{H}_{q^l}\rtimes \mathsf{C}_r$ such that $q,r$ are two odd primes, $l \geq 1$  and the action of $\mathsf{C}_r$ is faithful and simple. If $p=\min\{r,q\}$, then 
 $$\Pro(G)\leq \Pro(T_p)=\frac{1+\frac{p^2-1}{t(p)}}{p^2}=f_{n}(p).$$
\begin{proof}
We observe that all non-linear characters of $G$ have degree $r$. Therefore, we have that 
$$\Pro(T)=\frac{1+\frac{r^2-1}{q^l}}{r^2}.$$

Now, we have two possibilities for $p$: $p=r$ or $p=q$. Assume first that $r=p$. In this case, by definition of $t(p)$, we have that $q^l\geq t(p)$. Thus, we have that
 $$\Pro(T)= \frac{1+\frac{p^2-1}{q^l}}{p^2}\leq \frac{1+\frac{p^2-1}{t(p)}}{p^2}$$ 
 and the result holds in this case.

Assume now that $q=p$.  Since $r>p$ and $r,p$ are odd primes, we have that $r \geq p+2$. We claim that $l>2$. Assume by the contrary that $l=1$.  Then $r>p$ and $r$ divides $|\Aut(\mathsf{C}_p|=(p-1)$, which is impossible since $r\geq p+2$. Analogously, if $l=2$, then $r$ divides $|\Aut(\mathsf{C}_p\times \mathsf{C}_p)|=p(p+1)(p-1)^2$, which is impossible once again.  Thus, $l \geq 3$ and $ q^l =p^l \geq p^3$. Therefore, we have the following inequalities
$$\Pro(T)=\frac{1+\frac{r^2-1}{q^l}}{r^2}\leq \frac{1+\frac{(p+2)^2-1}{p^3}}{(p+2)^2}\leq \frac{1}{p^2}\leq \frac{1+\frac{p^2-1}{t(p)}}{p^2}.$$
Thus, the result also holds in this case.
\end{proof}
\end{lem}

Now, we are prepared to prove  Theorem A.

\begin{proof}[Proof of Theorem A]
Let $G$ be a counterexample of minimal order.  By Lemma \ref{red2}, we may assume  that  $G$ has the form $\mathsf{H}_{q^l} \rtimes \mathsf{C}_r$, where $r,q$ are two primes dividing $|G|$. Let $t=\min \{r,q\}\geq p$. Thus, applying Lemma \ref{prev2}, we deduce that 
$$\Pro(G)\leq f_{n}(t) \leq g_{n}(p)$$
and the result follows. 

It only remains to prove that the bound is sharp. Let $q\geq p$ be a  prime such that $g_{n}(p)=f_{n}(q)$ and let  $G:=\mathsf{C}_p\times T_q$. Then $G$ is non-nilpotent, $p$ is the smallest prime dividing $|G|$ and 
$$\Pro(G)=\Pro(\mathsf{C}_p)\Pro( T_q)=\Pro(T_q)=f_{n}(q)=g_{n}(p).$$
Thus, the bound is sharp.
\end{proof}

\subsection{Supersolvable odd groups}\label{supersolvable}


Our goal in this subsection is to prove Theorem B. As in the case of nilpotent groups, we will begin by reducing the problem to an special case. We say that $G$ is a just non-supersolvable group if $G$ is solvable but not supersolvable and every proper quotient of $G$ is supersolvable. We recall the Fitting subgroup of $G$ is defined as the largest normal and nilpotent subgroup, and we will denote it by $F(G)$. We have the following result.

\begin{thm}[Theorem 3.3 of \cite{RW}]\label{tar}
Let $G$ be a just non-supersolvable group. Then $G=F(G)\rtimes M$, where $M$ is a supersolvable group acting faithfully on $F(G)$ and $F(G)$ is the unique minimal normal subgroup of $G$. In particular, $F(G)$ is a non-cyclic elementary  abelian group and it is a faithful and simple $M$-module (equivalently, $G/F(G)$-module).
\end{thm}

In Theorem B we will only consider groups of odd order and hence the solvabilty of the groups involved will be guaranteed by the Feit-Thompson Theorem \cite{FT}. Thus, we will assume the solvability hypothesis in all preliminary results.

\begin{lem}\label{red1}
Let $G$ be solvable group which is not supersolvable group. Then $\Pro(G) \leq \Pro(T)$, where $T$ has the form $T=\mathsf{H}_{q^l}\rtimes \mathsf{C}_r$, where $q$ and $r$ are two primes dividing $|G|$, $l \geq 2$  and the action of $\mathsf{C}_r$ is faithful and simple.
\begin{proof}
Let $N$ be largest normal subgroup of $G$ such that $G/N$ is non-supersolvable. Then $G_0:=G/N$ is a just non-supersolvabe group and hence, by Theorem \ref{tar}, we have that $G_0=F(G_0)\rtimes (G_0/F(G_0))$, where $F(G_0)$ is $q$-elementary abelian for some prime $q$. Now, taking quotients and subgroups in $G_0$, we can obtain $T=\mathsf{H}_{q^l} \rtimes \mathsf{C}_r$, where the action of $\mathsf{C}_r$ is faithful and simple.  Thus, by Lemma \ref{basic}, we have the following inequalities
$$\Pro(G)\leq \Pro(G_0)\leq  \Pro(T)$$
and hence,  the result follows.
\end{proof}
\end{lem}

Now, we define $R_p:=\mathsf{H}_{r(p)}\rtimes \mathsf{C}_p$.  We observe that $R_p$ has the form described in Lemma \ref{red1} and that $\Pro(R_p)=f_{s}(p)$.

\begin{lem}\label{prev1}
Let $T=\mathsf{H}_{q^l} \rtimes \mathsf{C}_r$ such that $q,r$ are two odd primes, $l \geq 2$  and the action of $\mathsf{C}_r$ is faithful and simple. If $p=\min\{r,q\}$, then $$\Pro(T)\leq \Pro(R_p)=\frac{1+\frac{p^2-1}{r(p)}}{p^2}=f_{s}(p).$$
\begin{proof}
We have two possibilities for $p$: $p=r$ or $p=q$. If $q=p$, we can argue as in Lemma \ref{prev2} to prove that $l\geq 3$. Therefore, we can prove that 
$$\Pro(G)\leq \frac{1}{p^2}\leq  \Pro(R_p).$$
Finally, the case $r=p$ follows reasoning as in Lemma \ref{prev2}. 
\end{proof}
\end{lem}

Now, we finish this section by proving  Theorem B.

\begin{proof}[Proof of Theorem B]
Let $G$ be  not supersolvable of odd order. By Feit-Thompson Theorem \cite{FT}, we have that $G$ is a solvable group and hence, by Lemma \ref{red1}, we have that $\Pro(G) \leq \Pro(T)$, where $T$ has the form $T=\mathsf{H}_{q^l}\rtimes \mathsf{C}_r$ such that $q,r$ are two primes dividing $|G|$, $l \geq 2$  and the action of $\mathsf{C}_r$ is faithful and simple. Now, let $t=\min \{q,r\}\geq p$. Thus,  applying  Lemma \ref{prev1}, we have that 
$$\Pro(T)\leq f_{s}(t) \leq g_{s}(p)$$
and the result follows.

It only remains to prove that the bound is sharp. Let $q\geq p$ be prime such that $g_{s}(p)=f_{s}(q)$ and let $G:=\mathsf{C}_p\times R_q$. Then $G$ is non-supersolvable,  $p$ is the smallest prime dividing $|G|$ and 
$$\Pro(G)=\Pro(\mathsf{C}_p)\Pro( R_q)=\Pro(R_q)=f_{s}(q)=g_{s}(p).$$
Thus, the bound is sharp.
\end{proof}

\section{On the average character degree}\label{acd}

In this section we prove Theorems C and D. We start with the definition of the functions $h_n(p)$ and $h_s(p)$.
Let $p>2$ be a prime. We define $k(p)$ as the smallest positive integer such that the smallest  prime divisor of $k(p)$ is at least $p$ and $2k(p)+1=q^r$, where $q\geq p$ is a prime. In a similar way, we define $l(p)$ as  the smallest positive integer such that the smallest  prime divisors of $l(p)$ is at least $p$ and $2l(p)+1=q^r$, where $q\geq p$ is a prime and $r\geq 2$. In Section \ref{numbers}, we will prove that there exist integers, which satisfy the conditions defining $k(p)$ and $l(p)$ and hence both integers exist. With this notation we define
$$h_{n}(p)=\frac{3k(p)}{k(p)+2}$$
and
$$h_{s}(p)=\frac{3l(p)}{l(p)+2}.$$

We observe that $h_{n}(p),h_{s}(p)< 3$ for every prime $p$. Now, we prove Theorems C and D. We remark that our proofs are a refinement of the proof of \cite[Theorem D]{ILM}. As in case of the commuting probability, all groups considered here will have odd order and hence all of them will be solvable. During the proof of Theorems C and D we will use the Frattini subgroup. The Frattini of a group $G$ is defined as the intersection of all maximal subgroups of $G$ and it will be denoted by $\Phi(G)$.

\begin{proof}[Proof of Theorem C]
We observe that $p$ is a Sophie Germain prime if and only if $k(p)=p$ and hence, in this case, $h_{n}(p)=\frac{3p}{p+2}$. Thus, if $p$ is a Sophie Germain prime, then the result follows by Theorem \ref{original}. As a consequence, we may assume that $p\geq 7$.

Let $p\geq 7$ be a prime and let $G$ be a non-nilpotent group whose smaller prime divisor is $p$, $\acd(G)<h_{n}(p)$ and it has minimal order with these properties.  Let $A$ be a minimal normal subgroup of $G$. Since $G$ is solvable, we have that $A$ is an elementary abelian $l$-group for a prime $l$ (note that $l\geq p$) and since $G$ is non-nilptent we have that $A<G$.   Now, we have that $\acd(G)<h_{n}(p)<3<p$ and hence, by  Lemma \ref{dec}, we have that $\acd(G/A)\leq \acd(G)$ and, by minimality of $G$, we have that $G/A$ is nilpotent. 

If $A\leq Z(G)$, then $G/Z(G)$ is nilpotent and hence $G$ is nilpotent, which is impossible. Thus, we have that $1<[G,A]\leq A$ and $[G,A]\trianglelefteq G$. Since $A$ is  a minimal normal subgroup, we deduce that $[G,A]=A$. In addition, if $A\leq \Phi(G)$, then, applying \cite[Satz III.3.7]{Huppert}, we have that $G$ is nilpotent, which is impossible. Thus, there exists  a maximal subgroup $M$ of $G$, such that $A\not \leq M$, which implies that $G=MA$. We claim that $M\cap A=1$. Since $A\trianglelefteq G$, we know that $M\cap A\trianglelefteq M$ and since $A$ is abelian we know that $M\cap A\trianglelefteq A$. Thus,  $M\cap A \trianglelefteq MA= G$ and since $A\not \leq M$, the claim follows. 

 Thus, $G=A\rtimes M$ and hence, using Theorem \ref{thm32}, we can deduce that 
$$\frac{t(r+1)}{r+t}=f(t,r)\leq \acd(G)<h_{n}(p),$$
where $r$ is the number of orbits of the action of $G$ in $\Irr(A)\setminus \{1_{A}\}$ and $t$ is the size of one of these orbits. By Lemma \ref{reven}, we have that $r$ is even. Moreover, since   the action of $G$ on $A$  has no non-trivial fixed points, then applying \cite[Corollary 6.33]{Isaacscar}, we deduce that the action of $G$ in $\Irr(A)\setminus \{1_{A}\}$ has no orbits of size $1$. Therefore, $t>1$ and hence $t\geq p$.

Assume first that $r\geq 4$. Then $3 > h_{n}(p)> \acd(G)\geq f(p,4)\geq f(7,4)=\frac{35}{11}>3$, which is a contradiction. Thus, we may assume that $r=2$ and hence, applying  Lemma \ref{reven}, we have that both orbits have size $t$. Therefore, $|A|=2t+1$, where $|A|$ is the power of a prime at least  $p$, and the smallest prime divisor of $t$ is at least $p$. If $t$ is an integer satisfying this conditions, then $t\geq k(p)$. Thus, we have that
$$h_{n}(p)=\frac{3k(p)}{k(p)+2}\leq \frac{3t}{t+2}\leq \acd(G)<h_{n}(p),$$
which is again a contradiction.

Now, it only remains to prove that the bound is sharp. Let  $G=\mathsf{C}_p\times (\mathsf{H}_{2k(p)+1}\rtimes \mathsf{C}_{k(p)})$. Then $G$ is non-nilpotent, $p$ is the smallest prime dividing $|G|$ and we have that 
$$\acd(G)=\acd(\mathsf{C}_{p})\acd(\mathsf{H}_{2k(p)+1}\rtimes \mathsf{C}_{k(p)})=\acd(\mathsf{H}_{2k(p)+1}\rtimes \mathsf{C}_{k(p)})=\frac{3k(p)}{k(p)+2}=h_{n}(p).$$
Thus, the bound is best possible.
\end{proof}

\begin{proof}[Proof of Theorem D]
Let $p>3$ be a prime and let $G$ be a non-supersolvable group whose smaller prime divisor is $p$, $\acd(G)<h_{s}(p)$ and it has minimal order with these properties. Let $A$ be a minimal normal subgroup of $G$. As in the proof of Theorem C, we have that $A$ is an elementary abelian $l$-group for a prime $l$  and that $A<G$. Again, we have that $\acd(G)<h_{s}(p)<3<p$ and hence, by  Lemma \ref{dec}, we have that $\acd(G/A)\leq \acd(G)$ and by minimality of $G$, we have that $G/A$ is supersolvable.  Since $G$ is non-supersolvable, we have that $A$ cannot be cyclic and  $A \not \leq Z(G)$, which implies $[G,A]=A$. On other side, if $A\leq \Phi(G)$, then $G/\Phi(G)$ is supersolvable and hence, applying \cite[Satz VI.8.6(a)]{Huppert}, we have that $G$ is supersolvable, which is a contradiction.  Thus, reasoning as in Theorem C,  we can prove that $G$ splits over $A$ and hence, applying Theorem \ref{thm32}, we have that 
$$f(t,r)\leq \acd(G)<h_{s}(p),$$
where $r$ is the number of orbits of the action of $G$ in $\Irr(A)\setminus \{1_{A}\}$ and $t$ is the size of one of these orbits.  Reasoning as in proof of Theorem C, we can deduce that $t\geq p$, and $r\geq2$ is even.

Assume first that $r \geq 4$. If $t \geq 7$, then  
$$3<\frac{35}{11}=f(7,4)\leq f(t,r)<h_{s}(p)< 3,$$
which is a contradiction. Thus, we may assume that $t<7$, which forces $p=5$ and $t=5$.

 We claim that, in this case, the size of all  orbits of the action in $\Irr(A)\setminus \{1_A\}$ is $5$.  Since $t=5$ is the size of an orbit, we have that there exists $\lambda\in \Irr(A)$ such that the size of the $G$-orbit containing  $\lambda$ is $5$. Let $T$ be the stabilizer of $\lambda$ in $G$. Then we have that $|G:T|=5$ and since $5$ is the smallest prime dividing $|G|$, we deduce that $T$ is normal in $G$. Therefore, $[A,T]\trianglelefteq G$ and since $[A,T]<A$, then we have that $[A,T]=1$. As a consequence, $T$ is contained in the stabilizer of all characters in $\Irr(A)$, which implies that the size of any orbit in $ \Irr(A)\setminus \{1_A\}$ is at most $5$ and the claim follows.  

Now, if $r \geq 6$, then 
$$3<\frac{35}{11}=f(5,6)\leq f(t,r)<h_{s}(p)< 3,$$
which is a contradiction. Thus, we have that $r=4$ and hence $|A|=4\cdot5+1=21$, which is impossible since $21$ is not a prime power.

Finally, assume  that $r=2$. By Lemma \ref{reven}, we have that  the two orbits have the same size $t$. Thus, $|A|=2t+1$, where $|A|$ is not prime, $|A|$ is the power of a prime at least $p$, and the smallest prime divisor of $t$ is at least $p$. If $t$ is an integer satisfying these conditions, then $t\geq l(p)$, by definition. Therefore, we have that
$$h_{s}(p)=\frac{3l(p)}{l(p)+2}\leq \frac{3t}{t+2}\leq \acd(G)<h_{s}(p),$$
which is a contradiction.

Now, it only remains to prove that the bound is sharp. Let  $G=\mathsf{C}_p\times (\mathsf{H}_{2l(p)+1}\rtimes \mathsf{C}_{l(p)})$. Then $G$ is non-supersolvable, $p$ is the smallest prime dividing $|G|$ and we have that 
$$\acd(G)=\acd(\mathsf{C}_{p})\acd(\mathsf{H}_{2l(p)+1}\rtimes \mathsf{C}_{l(p)})=\acd(\mathsf{H}_{2l(p)+1}\rtimes \mathsf{C}_{l(p)})=\frac{3l(p)}{l(p)+2}=h_{s}(p).$$
Thus, the bound cannot be improved.
\end{proof}

\section{Number-theoretical questions}\label{numbers}

We  begin by proving that numbers $k(p)$ and $l(p)$ defined in  Section \ref{acd} do indeed exist. Then, we will discuss some properties of  $t(p)$ and $r(p)$.

Before beginning, we define the recall the definition of the  primordial of a number. Given $n\in \mathbb{N}$ we define the primordial of $n$ as the product of all primes smaller than $n+1$, that is
$$n\#=\prod_{q\leq n,q\text{ is prime}}q.$$
It is easy to see that if $n$ is a positive integer and $p$ is a prime, then  the smallest prime divisor of $n$ is at least $p$ if and only if $(n,(p-1)\#)=1$.  Now, we can prove the existence of $l(p)$. We observe that if  there exists an integer satisfying the condition of $l(p)$, then this integer  satisfies the conditions of $k(p)$ and hence we will prove the existence of both numbers at the same time. The following result is due to  Bryce Kerr.

\begin{thm}[Kerr]
Let $p$ be a prime. Then there exists  $l \in \N$ satisfying the following two conditions.
\begin{itemize}
\item[(i)]  The smallest prime divisor of $l$ is at least $p$.

\item[(ii)] $2l+1=q^f$ where $q\geq p$ is a prime and $f\geq 2$.
\end{itemize}
\begin{proof}
Let $r$ be a prime. We observe that 
$$(r-1)(4(r-1)^2+6(r-1)+3) \equiv -1 \pmod{r}$$
and that
$$2(r-1)+1 \equiv -1 \pmod{r}.$$

Since $(p-1)\#=\prod_{r<p,r \text{ is prime}}r$, applying the Chinese Remainder Theorem, we have that there exists $a\in \N$ such that $n \equiv a \pmod{(p-1)\#}$ if and only if $n \equiv -1 \pmod{r}$ for every $r<p$, prime.

Now, we claim that there exists $n$ such that $n \equiv a \pmod{(p-1)\#}$ and $2n+1=q$, where $q$ is a prime larger than $p$. This is equivalent to showing that there exist a prime $q>p$, of the form
$$q=k 2 (p-1)\#+(2a+1)$$
for some $k \in \N$. Now, we have that $2a+1 \equiv 2(r-1)+1\equiv -1 \pmod{r}$ for every $r<p$, prime. Therefore, $(2(p-1)\#,2a+1)=1$ and hence the prime $q$ exists by Dirichlet's Theorem.

Now, let $n$ as in the previous claim and let $l=n(4n^2+6n+3)$. By the choice of $n$, we have that
$$q^3=(2n+1)^3=2n(4n^2+6n+3)+1=2l+1.$$

Thus, $l$ satisfies condition (ii) of the statement. Now, since $n \equiv a \pmod{(p-1)\#}$, we have that for every $r<p$ prime 
$$l\equiv n(4n^2+6n+3)\equiv (r-1)(4(r-1)^2+6(r-1)+3)\equiv -1 \pmod{r}.$$

This implies that $(l,(p-1)\#)=1$ and hence the smallest prime divisor of $l$ is at least $p$. Thus, $l$ satisfies condition  (i) of the statement and the result follows.
\end{proof}
\end{thm}

\begin{cor}[Kerr]
Let $p$ be a prime. Then there exists $k \in \N$ satisfying the following two conditions.

\begin{itemize}

\item[(i)] The smallest prime divisor of $k$ is at least $p$. 

\item[(ii)] $2k+1=q^f$ with $q\geq p$ and $f\geq 1$.
\end{itemize}
\end{cor}

Now, we discuss about $t(p)$ and $r(p)$. First, we introduce some notation. Given $a,d\in \mathbb{N}$ with $0<a<d$ and $(a,d)=1$ we define $P(a,d)$ as the smallest prime in the arithmetic progression $\{a+n d\}_{n \in \mathbb{N}}$ (such a prime exists by Dirichlet's Theorem).  Computations in GAP \cite{gap}  have shown that if $p\leq 1000$, then $t(p)=P(1,p)$ and  $r(p)=P(p-1,p)^2$. We conjecture that this holds for all primes.

\begin{con}\label{conAB}
Let $p$ be a prime. Then the following hold
\begin{itemize}

\item [a)] $t(p)=P(1,p)$.

\item [b)] $r(p)=P(p-1,p)^2$.
\end{itemize}
\end{con}

Conjecture \ref{conAB} depends on strong properties of the distribution of prime numbers. According to B. Kerr and T. Trudgian, it seems that Conjecture \ref{conAB} should be true, but it may be out of reach with current knowledge on number theory.  Linnik \cite{L} proved that there exist two positive constants $c,L$ such that $P(a,d)\leq c \cdot d^{L}$ for every pairs of numbers $a,d$. Many efforts have been done to give explicit constants $c$ and $L$. The best bound known for $L$ is due to Xylouris \cite{X} and it is $L=5$. However, it is  conjectured that $P(a,d)\leq d^2$ (see for example \cite{H}). It is not hard to see that if this conjecture holds, then part a) of Conjecture \ref{conAB} also holds.  It is worth to say that the bound $P(a,d)\leq d^2$ is not known even assuming the Riemann Generalized Hypothesis. We close this paper with a few more comments about conjecture \ref{conAB}.

\begin{note}
Fix a prime $p$. If we assume that $P(1,p),P(p-1,p)\leq p^k$ and $p$ a Zsigmondy prime for $ \langle q,r \rangle$ with $q>p$ and  $r\geq k$, then $q^r\geq p^k\geq  t(p)$ and if $r\geq 2k$, then $q^r\geq r(p)$. Thus, we have the following.
\begin{itemize}
\item[a)] Using Xylouris's bound \cite{X}, we have that $P(1,p),P(p-1,p)\leq p^5$. If $\langle q,r\rangle$ is a counterexample for part a) of Conjecture \ref{conAB}, then $r\in \{2,3,4\}$ and if it is a counterexample for part b), then $r \in \{3,4,5,6,7,8,9\}$.

\item[b)] Assuming the Generalized Riemann Hypothesis, we can use the result of   Lazmouri, Li and Soundararajan \cite{LLS}, to deduce that $P(1,p),P(p-1,p)\leq (p-1)^2(\log(p))^2<p^3$. Thus, if $\langle q,r\rangle$ is a counterexample for part a) of Conjecture \ref{conAB}, then $r=2$ and if it is a counterexample for part b) then $r \in \{3,4,5\}$. Note that the bound of \cite{LLS} is not enough to prove part a) of Conjecture \ref{conAB}. For $p=7$ we have that $r(7)=169\in \T(7,1)$, but $169<(7-1)^2(\log(7))^2$.
\end{itemize}

\end{note}

\bigskip

\centerline{\bf Acknowledgement}

\bigskip

The author would like to thank  Bryce Kerr and Timothy Trudgian for their help  in the  number-theoretical questions of this paper. The author would also like to thank Alexander Moretó for suggesting me this subject, and for helpful conversations.


\end{document}